\documentclass[12pt]{article}
\usepackage{amssymb} 

\def\i{\iota}

\def\o{\omega}

\def\u{\upsilon}

\def\S{\Sigma}

\chardef\tempcat=\the\catcode`\@
\catcode`\@=11
\def\cyracc{\def\u##1{\if \i##1\accent"24 i
    \else \accent"24 ##1\fi }}
\newfam\cyrfam

\DeclareFontFamily{OT1}{msb}{}{}
\DeclareFontShape{OT1}{msb}{m}{n}
 {  <5> <6> <7> <8> <9> <10> gen * msbm
      <10.95><12><14.4><17.28><20.74><24.88>msbm10}{}
\DeclareMathAlphabet{\bubble}{OT1}{msb}{m}{n}

\def\bR{{\bubble R}}

\def\bC{{\bubble C}}

\def\bG{{\bubble G}}
\def\bU{{\bubble U}}

\newfont{\goth}{eufm10 scaled \magstep1}

\def\gg{\mbox{\goth g}}
\def\gh{{\mbox{\goth h}}}

\def\gso{\mbox{\goth so}}

\def\gsp{\mbox{\goth sp}}

\def\SpV{{\mathrm{Sp({\it V})}}}

\def\SOV{{\mathrm{SO({\it V})}}}

\def\GL#1{{\mathrm{GL(#1)}}}
\def\GLV{{\mathrm{GL({\it V})}}}
\def\SLV{{\mathrm{SL({\it V})}}}
\def\Aff#1{{\mathrm{Aff(#1)}}}
\def\AffV{{\mathrm{Aff({\it V})}}}
\def\AffGV{{\mathrm{Aff_{\it G}({\it V})}}}
\newfont{\mcal}{eusm10 scaled \magstep1}

\def\square{\kern1pt\vbox
            {\hrule height 0.6pt\hbox{\vrule width 0.6pt\hskip 3pt
 \vbox{\vskip 6pt}\hskip 3pt\vrule width 0.6pt}\hrule height 0.6pt}\kern1pt}

\def\ra{\rightarrow}

\def\ker{\mathrm{ker\;}}

\newtheorem{Th}{Theorem}
\newtheorem{Prop}{Proposition} 
\newtheorem{Cor}{Corollary} 
\newtheorem{Lem}{Lemma}
\newtheorem{Def}{Definition} 

\def\bt{\begin{Th}}
\def\et{\end{Th}}
\def\bp{\begin{Prop}}
\def\ep{\end{Prop}}
\def\bc{\begin{Cor}}
\def\ec{\end{Cor}}
\def\bl{\begin{Lem}}
\def\el{\end{Lem}}
\def\bd{\begin{Def}}
\def\ed{\end{Def}}

\def\pf{\noindent{\it Proof. }}
\def\pfof#1{\noindent{\it Proof of #1. }}
\def\qed{\hspace{2ex} \hfill\square \par \medskip}  
\def\n{\nabla}
 \def\ot{\otimes}

\def\be{\begin{equation}}
\def\ee{\end{equation}}

\def\arr{\begin{array}{rlll}}
\def\ea{\end{array}}
\def\bea{\begin{eqnarray}}
\def\eea{\end{eqnarray}}  
\def\bean{\begin{eqnarray*}}
\def\eean{\end{eqnarray*}}  

\catcode`@=11
\@addtoreset{equation}{section}
\catcode`@=12   

\def\3s{3-Sasakian manifold}

\begin{document}

\title{Abelian simply transitive affine groups\\ of symplectic type} 
\date{December 7, 2001}

\author{Oliver Baues \thanks{e-mail: oliver@math.ethz.ch} \\
Departement Mathematik\\ 
ETH-Zentrum \\ 
R\"amistrasse 101 \\
CH-8092 Z\" urich 
\and
Vicente Cort\' es \thanks{e-mail: vicente@math.uni-bonn.de \newline
This work was supported by FIM (ETH Z\"urich) and  
MPI f\"ur Mathematik (Bonn).} \\
Mathematisches Institut \\
Universit\"at Bonn\\
Beringstra\ss e 1 \\
D-53115 Bonn
}
\maketitle
\begin{abstract} \noindent  
The set of all Abelian simply transitive subgroups
of the affine group naturally corresponds to the set 
of real solutions of a system of algebraic equations. We
classify all simply transitive subgroups of 
the symplectic affine group by constructing 
a model space for the corresponding variety of solutions.
\noindent 
Similarly,  we classify the complete global 
model spaces for flat special K\"ahler manifolds with
a constant cubic form. 
 

\end{abstract} 
\section{Introduction} 
Affine groups acting simply transitively on real affine space have been much 
studied in the literature.  (See e.g.\ the seminal paper \cite{A}, or \cite{S} for
a different point of view.)  In this paper we study
the case of Abelian groups which act by 
symplectic transformations. Our interest is to describe
the properties of the {\em algebraic variety\/} of {\em simply
transitive Abelian affine groups of symplectic type}.
By doing so, we solve the classification problem for 
these groups up to a linear action on a well understood algebraic space.  
To put this into perspective we mention that in general Abelian simply
transitive groups of affine motions are classified only in low dimensions.
(See \cite{DO} for some recent contribution to this problem.)
On the other hand, it is well known that the group of translations is
the only simply transitive Abelian group with the property that 
its linear part preserves a nondegenerate symmetric bilinear form.
We mention further that (as follows for example from 
\cite[Theorem 6.11,7.1]{FGH})
the classification of
Abelian simply transitive affine groups of symplectic type 
essentially implies the classification of flat affine tori with
a parallel symplectic form up to symplectic affine diffeomorphisms. 


Our particular motivation to consider Abelian simply
transitive symplectic affine groups stems from {\em special K\" ahler geometry}. 
Special K\" ahler geometry is  a particular type of geometry which arises
in certain supersymmetric field theories, see e.g.\ 
\cite{F} and the literature cited therein.  
A special  K\" ahler manifold 
is a  (pseudo-) K\" ahler manifold $(M,g)$ which has as additional 
geometric datum
a torsionfree flat connection $\nabla$ which satisfies certain compatibility
conditions with respect to $g$.   
A special K\"ahler manifold is called 
flat if the Levi-Civita connection $D$ of the special K\"ahler metric $g$
is flat. In this paper, we are 
concerned with the existence and construction of {\em flat\/} special 
K\"ahler manifolds. 

Clearly, if a K\"ahler manifold is flat, it is locally modelled 
on the vector space $\bC^n$ endowed with a Hermitian inner product, 
and therefore it is locally trivial as a K\"ahler manifold.
We say that a special K\"ahler manifold is locally trivial 
if as a {\em special\/} K\"ahler manifold it is locally 
equivalent to a Hermitian vector space. This is the case 
precisely if $\nabla=D$. 
It may seem, at first sight, a bit surprising that there do exist
flat special K\"ahler manifolds which are not locally 
trivial as special K\"ahler manifolds. 
In fact, we show that if the special K\"ahler metric is
definite, then flatness implies local triviality. In particular, 
any simply connected and complete flat special K\"ahler manifold with a
positive definite metric is equivalent to a positive definite
Hermitian vector space. 

It is well known that special K\" ahler manifolds arise 
locally from holomorphic potentials.
As it turns out, flat special K\" ahler manifolds correspond 
to those holomorphic potentials which satisfy
a specific algebraic constraint on their third derivatives. 
(A priori, flatness is a constraint on the 
fourth order jet of the potential.)   
We may interpret this constraint as an 
associativity condition on an associated bilinear product
on the tangent spaces of a special K\"ahler manifold. 
This shows, by the way,    
a close structural similarity between {\em flat 
special K\" ahler manifolds} and {\em Frobenius manifolds}. 
In the non-definite case, the algebraic constraint 
on cubic tensors admits a non-trivial variety of solutions
which we describe as certain bundles over Grassmannians. 
This enables us to locally classify flat special 
K\"ahler manifolds with constant cubic form. In particular, 
the points of the above solution variety correspond to the 
complete simply connected global models of such manifolds.
These spaces are complete with respect to both the metric connection $D$,
and the affine connection $\nabla$. All of them 
arise from certain Abelian simply transitive affine groups of symplectic 
type. 

\tableofcontents 

\subsection{Statement of main results and outline of the paper} 
\paragraph{The variety of Abelian simply transitive groups}
The set of simply transitive
Abelian groups which act affinely on a vector space $V$ 
with linear part in a prescribed linear Lie group $G \subset \GLV$ 
constitutes  in a natural way a real affine variety. 
To explain this we show that this set corresponds 
to the solutions  $C(\gg)$ of a system
of linear and quadratic equations in the first prolongation 
\[\gg^{(1)} := V^*\otimes \gg \cap S^2(V^*) \otimes V\] 
of the Lie algebra $\gg = Lie\, G$. This is done in section \ref{av}. 
 
Our guiding principle is that 
we want not only to write down the equations defining this variety, namely, 
\[ C(\gg ) = \{ S \in \gg^{(1)} \mid {\rm tr}\,  S_X = 0 \quad \mbox{and} 
\quad [S_X,S_Y] =0 \quad \mbox{for all} \quad X,Y \in V\} \]  
but want also to {\em solve} them by some explicit construction.
This is carried out in section~\ref{ospSec}. 

We need to introduce some notation to explain our result for  $\gg= \gsp(V)$. 
Let $V$ be a vector space with a nondegenerate alternating product $\omega$,
$\dim V = 2n$.  A subspace $U$ of $V$ is called {\em Lagrangian} if $U$ is maximally isotropic
with respect to $\omega$. For $S \in \gg^{(1)}$ we define the {\em support of $S$}
by \[ \S_S := {\rm span} \{ S_XY \mid  X,Y \in V\} \; .\] 
The variety $C(\gsp(V))$ is naturally stratified by the dimension of the support. 
The symplectic form $\omega$ allows the identification of $S^3 V$
with a subspace of $V^* \otimes V^* \otimes V$ and a corresponding
notion of support for elements of $S^3 V$. For a subspace $W \subset V$
we define 
\[    S^3W_{reg} := \{ S\in S^3W\mid \S_S = W\} \; . \]  
The first result is (compare Theorem \ref{1stmainThm}):
\bt Let  $U\subset V$ be  a Lagrangian subspace.
There is a  one-to-one correspondence between 
Abelian simply transitive affine groups   
of symplectic type up to conjugation in the affine symplectic group and 
elements of $S^3U$ up to $\SpV$-equivalence (two elements of $S^3U$ 
are $\SpV$-equivalent if they lie on the same $\SpV$-orbit in 
$S^3V$).   
\et

Let ${\cal U}_k \ra G_k^0(V)$ denote the universal (tautological) 
vector bundle over the  Grassmannian $G_k^0(V)$ of $k$-dimensional isotropic subspaces, 
and $S^3{\cal U}_k$ the third symmetric power of this bundle. We prove then
(compare Theorem \ref{stratiThm}):  

\bt  
The $k$-th stratum $C(\gsp(V))_k$ of the variety $C(\gsp(V))$ naturally identifies
with the Zariski open subbundle
\be 
\bigcup_{W \in G_k^0(V)} S^3W_{reg}  \ee
of the vector bundle $S^3{\cal U}_k$. 
\et

Now the classification for Abelian simply transitive affine groups   
of symplectic type reads as follows 
(compare Theorem \ref{stratiThm}, Corollary \ref{translCor}):
\bc
There is a natural one-to-one correspondence between  
Abelian simply transitive affine groups   
of symplectic type with (maximal) translation subgroup of dimension $2n-k$ up
to conjugation, and homogeneous cubic polynomials in 
$S^3W_{reg}$ up to the linear action of $\GL W$, where $W$ is a
$k$-dimensional isotropic subspace of $V$. 
\ec
In particular, every such group contains
an  $n$-dimensional subgroup of translations. \newline

\paragraph{Special K\"ahler manifolds}
Let us turn now to the second part of our article which starts in 
section \ref{fskm}. 
A special K\"ahler manifold
is a (possibly indefinite) K\"ahler manifold $(M,J,g)$ endowed with a 
flat torsionfree connection $\nabla$ such that $\nabla J$ is symmetric
and $\n \o = 0$, where $\o = g (\cdot , J \cdot )$ is the symplectic 
(K\"ahler) form \cite{F}. A special K\"ahler manifold is 
called flat if the Levi-Civita connection $D$ is flat. The simplest 
example of such a manifold is obtained by taking
a (pseudo-) Hermitian vector space $\bC^n$, so that
the flat Levi-Civita connection $D$ and the symplectic
connection $\nabla$ coincide. Henceforth a special K\"ahler manifold
is called trivial if $D = \nabla$. 
We show:

\bt Any flat special K\"ahler manifold with a definite metric
is trivial. 
\et
This result is obtained by considering the algebraic
constraints on the difference tensor $S= D - \nabla$ 
of a flat special K\"ahler manifold.  Let 
$(V = \bC^n,J,g)$ be the standard
\mbox{(pseudo-)} Hermitian vector space of complex signature $(p,q)$, $p+q = n$.  
We show that the
tensor field $S$ takes values in a subvariety 
$C_J(\gsp(V))$ of the cone $C(\gsp(V))$. 
In section \ref{vc} we obtain 
a description of the variety $C_J(\gsp(V))$ in a way
analoguous to Theorem 1 and Theorem 2 above. These
results allow for the local classification of 
special  K\"ahler manifolds in section \ref{lc} as follows: 

\bt Let $f$ be a holomorphic function defined on an open 
subset $M \subset V$.
Assume that the pointwise support of the cubic tensor field 
$\partial^3f : M \ra  S^{3,0}V^*$ defined by the holomorphic  
third partial derivatives of $f$ is isotropic and put $S_f := \partial^3f + 
\overline{\partial^3f} \in S^3V^*$. Then $M_f := 
(M,J,g, \n := D + S_f)$
is a flat special K\"ahler manifold and any flat special K\"ahler manifold
arises locally in this way. 
\et 

As a special case of the theorem we can describe the 
class of flat special K\"ahler manifolds with constant 
cubic form, i.e., those manifolds which satisfy $DS = 0$. 
They admit a model with an Abelian simply 
transitive group of automorphisms. 

\bt 
Let $W\subset V$ be a complex isotropic subspace, and $f: V \ra \bC$ 
a holomorphic cubic polynomial which identifies with an element of $S^{0,3}W$ 
under the canonical identification $S^{0,3}V^* \cong S^{3,0}V$. 
Then the manifold 
$V_f$, defined in the previous theorem, is a flat special 
K\"ahler manifold with constant cubic form and complete connections $D$ and 
$\n$. Moreover, any 
flat special K\"ahler manifold with 
constant cubic form is locally equivalent 
to an open subset of some space $V_f$.
\et 
The result may be refined a bit. 
We assume that $(V,J,g)$ admits a complex 
Lagrangian subspace $U$.
Then any flat special K\"ahler manifold with 
constant cubic form and without trivial factor
may be obtained by some $S \in S^{0,3}U_{reg}$,
and we obtain a bijection between the orbits of the group 
${\rm GL}_{\bC}(U)$ on  
$S^{0,3}U_{reg}$ and equivalence classes of 
germs of flat  special K\"ahler manifolds with constant 
cubic form and without trivial factor. (See Theorem \ref{5thmainThm}.) 
Any such manifold has complex signature $(m,m)$, $m = \dim_{\bC} U$. 

\section{Abelian simply transitive affine groups} \label{As}
Let $V$ be a real vector space and $G \subset \GLV$ any Lie   
subgroup. We denote by $\AffGV \subset \AffV$
the group of affine transformations with linear part in $G$. 
Connected Lie subgroups of $\AffGV$ are called 
{\em affine groups of type} $G$. 
Lie subalgebras of ${\rm aff}_{\gg}(V) = Lie\, \AffGV$ are called 
{\em affine Lie algebras of type} $\gg = Lie\, G$.
We are interested in Abelian affine groups $H$ of type $G$ acting {\em 
simply transitively} on $V$. This means that the orbit map 
$\varphi : H \ni h \mapsto 
h0 \in V$ is a diffeomorphism. More generally, we will consider {\em 
almost simply transitive} groups, i.e.\ groups for which $\varphi$ is an 
immersion. The corresponding Lie subalgebras $\gh = Lie\, H \subset 
{\rm aff}_{\gg}(V)$ are almost  simply transitive.  
This means that the linear map $\phi : 
\gh \ni X \mapsto X0 \in V$ is an isomorphism.
In fact, this map is the differential of the orbit map $\varphi$, which is 
an immersion since $H$ is almost simply transitive. 
The Lie algebras of simply transitive groups satisfy a stronger condition. 
They are simply transitive, which means that for all $v\in V$ the linear map 
$\gh \ni X \mapsto Xv \in T_vV \cong V$ is an isomorphism. 
The next proposition reduces the study of Abelian (almost) 
simply transitive affine 
groups to that of Abelian (almost) simply transitive affine Lie algebras. 

\bp The $Lie$-functor from the category of Lie groups to the category of 
Lie algebras 
induces a bijection between the set of 
Abelian simply transitive affine groups of type $G$ 
and the set of Abelian simply transitive affine 
Lie algebras of type $\gg = Lie\, G$. The same is true for almost
simply transitive affine groups and Lie algebras.   
\ep

\pf We have seen that the Lie algebra of an Abelian (almost) 
simply transitive group 
of type $G$ is an Abelian (almost) 
simply transitive Lie algebra of type $\gg$. It is clear that 
the Abelian affine group $H$ of type $G$ generated 
by an Abelian almost simply transitive Lie algebra $\gh$ 
of type $\gg$ is almost simply transitive. 
It remains to show that the  group generated 
by a simply transitive Lie algebra $\gh$ 
is simply transitive. {}From the fact that $\gh$ is simply 
transitive it follows that all orbits of $H$ are open. $V$ being connected,
this implies that $H$ is transitive. So we have a diffeomorphism 
$V \stackrel{\sim}{\ra} H/H_v$, where $H_v$ is the 
stabilizer of a point $v\in V$. Moreover $H_v$ is discrete and $H \ra H/H_v 
\cong V$ is a covering of the simply connected manifold $V$. This implies 
that $H_v$ is trivial and hence that $H$ is simply transitive. 
\qed   

Let us denote by $A'(G)$ (respectively $A'(\gg )$) the set of 
Abelian almost simply transitive affine groups of type $G$ (respectively 
the set of Abelian almost simply transitive affine 
Lie algebras of type $\gg$). The subsets consisting of simply transitive 
groups and Lie algebras are denoted by $A(G)$ and $A(\gg )$, respectively.
Note that by the previous proposition we can identify $A(G) = A(\gg )$
and $A'(G) = A'(\gg )$. 
  
\subsection{The variety of simply transitive groups} \label{av}
Let $\iota : \gh  \hookrightarrow  {\rm aff}_{\gg}(V)$ be 
the canonical inclusion map associated to a 
Lie algebra $\gh \in A'(\gg)$. It gives rise  
to a injective linear map  $\rho = \rho_{\gh} : V \longrightarrow {\rm aff}_{\gg}(V)$ by 
\be \rho := \iota \circ \phi^{-1} \, \ee
where $\phi$ is the differential of the 
orbit map $\varphi: H \ni h \mapsto h0 \in V$ at the identity. 
We remark that, since $\gh$ is assumed to be an Abelian Lie algebra, 
the linear map $\rho: V \longrightarrow {\rm aff}_{\gg}(V)$ is in fact a
homomorphism of Lie Algebras.
 The affine Lie algebra ${\rm aff}_{\gg}(V) = \gg + V$ is the semidirect sum 
of the linear 
Lie algebra $\gg$ and the ideal $V$ of infinitesimal translations. 
We denote an element of this semidirect sum as a pair $(S,t)$, where
$S \in \gg$ is the linear part and $t \in V$ the translational part. 
In particular we can write the monomorphism $\rho = \rho_{\gh}$ in the form 
\be \rho (X) = (S_X, t_X)\, , \quad X\in V \, , \ee
where $S: V \ra \gg$ and $t :  V \ra V$ are linear maps. 

\bp \label{i-iiProp} 
Let $\rho = (S,t) : V \longrightarrow {\rm aff}_{\gg}(V) = \gg + V$ be 
the monomorphism associated to an Abelian almost 
simply transitive affine Lie algebra
of type $\gg$. Then $t = {id}_V$ and $S\in V^*\otimes \gg$ satisfies: 
\begin{itemize}
\item[(i)] $S_XY = S_YX$,  for all $X,Y \in V$. 
\item[(ii)] $[S_X,S_Y] = 0$,   for all $X,Y \in V$.
\end{itemize}
Conversely, any $S\in V^*\otimes \gg$ satisfying (i) and (ii) 
defines a monomorphism $\rho = (S,{id}_V): V \longrightarrow  {\rm aff}_{\gg}(V)$
onto an Abelian almost simply transitive affine Lie algebra 
$\gh = \rho (V)$ of type $\gg$.  $\gh$  is simply transitive if and only if   
\begin{itemize}
\item[(iii)]  ${\rm tr}\,  S_X = 0$,  for all $X \in V$.  
\end{itemize}
\ep

\pf Let us first check that $t = {id}_V$. This follows from 
\[t_X = \rho (X)0 = \iota (\phi^{-1}(X))0 = \phi^{-1}(X)0 = 
\phi (\phi^{-1}(X)) = X\, .\]  
Equations (i-ii) express that $\rho : V \ra 
{\rm aff}_{\gg}(V)$ is a homomorphism of Lie algebras.  

Now we show that (iii) is equivalent to the image of  $\rho$ being a  
simply transitive Lie algebra.  Let us first assume that $\gh = \rho (V)$  is a   
simply transitive Lie algebra. Let $H$ be the corresponding simply
transitive subgroup of $\Aff V$. For $h \in H$ we consider the
Jordan-decomposition $h = h_s h_u$ inside the real linear algebraic group
$\Aff V$.  Since $H$ is Abelian, $T = \{ h_s \mid h \in H \}$ is 
a group of semisimple operators which centralizes $H$. By semi-simplicity 
$T$ has a fixed point on $V$. But since $T$ centralizes the transitive
group $H$ it must be trivial. Therefore $H$ is unipotent, and in particular
the linear parts of the elements of $H$ are unipotent. This
implies (iii). 

Conversely, we prove that (i)-(iii) implies that $\gh$  is simply transitive. We have to 
show that the linear map
$V\ni X \mapsto S_Xv + X \in V$ is an isomorphism for all $v \in V$. Therefore it is sufficient to 
prove that the map $X \mapsto S_Xv = S_vX$ is nilpotent, for all $v \in V$.
 Let us show that ${\rm tr}\, S_v^k = 0$ for all
$k = 1,2, \ldots \, $.  Because of (i) and (ii) we have $S_v^k  = S_w$, 
where $w :=  S_v^{k-1}v$.  Now ${\rm tr}\,  S_v^k = 0$ follows from (iii). This shows that
$S_v$ is nilpotent, for all $v \in V$. 
 \qed

We recall \cite{Kob} that  the first 
prolongation of a Lie algebra $\gg \subset V^* \otimes V$ is defined as 
\[\gg^{(1)} := V^*\otimes \gg \cap S^2(V^*) \otimes V\, .\]   
Note that the equation (i) says that $S$ is an element of $\gg^{(1)}$. 
We denote by $C'(\gg ) \subset \gg^{(1)}$
the cone defined by the system (ii) 
of homogeneous quadratic equations. It is an
affine real algebraic variety. The condition (iii) defines a Zariski closed 
subset $C(\gg)  \subset C'(\gg)$ . In fact, $C(\gg)  = C'(\gg_0)  =  
C'(\gg) \cap \gg_0^{(1)}$, where 
$\gg_0 = \{ A\in \gg \mid {\rm tr}\,  A = 0\}$.

\bt \label{generalThm} 
The correspondence $H \mapsto S$, where $S \in \gg^{(1)}$ is the linear part of 
$\rho_{\gh}$,  induces bijections $A'(G) \stackrel{\sim}{\ra} C'(\gg )$
and $A(G) \stackrel{\sim}{\ra} C(\gg )$. 
This defines on the set $A'(G)$ of Abelian almost simply transitive affine 
Lie groups of type $G$ the structure of an affine cone over
a quadratic projective real algebraic variety. The subset 
$A(G) \subset A'(G)$ consisting of simply transitive groups 
is a closed subvariety defined by a system of linear equations.    
Under the above identification,  
the action of $g\in G$ on $A'(G)$ by conjugation 
corresponds to the linear 
action $\gg^{(1)} \ni S \mapsto g\cdot S \in \gg^{(1)}$, 
where $(g\cdot S)_X := gS_{g^{-1}X}g^{-1}$. 
\et  

\pf This follows essentially from Proposition \ref{i-iiProp}. We check
the formula for the $G$-action on $C'(\gg ) \subset \gg^{(1)}$. {}From the
definition of $g\cdot S$ it follows that the homomorphism 
$((g \cdot S),{id}_V)$ has the same image as 
\[ X \mapsto (gS_Xg^{-1}, gX) = g\rho_{\gh}(X)g^{-1} \, ,\]
namely the conjugated Lie subalgebra  
${\rm Ad}_g(\gh ) = g\gh g^{-1} \subset {\rm aff}_{\gg}(V)$. This shows that 
$((g \cdot S),{id}_V) = \rho_{{\rm Ad}_g(\gh )}$ and hence that
$S^{gHg^{-1}} = g\cdot S^H$. 
\qed

Now we specialize to the unimodular case, i.e.\ we assume that  
$G \subset \SLV$.

\bc \label{fundCor} An 
Abelian almost simply transitive affine group $H$ of 
of type $\SLV$ is simply transitive. In particular, $A(G) = A'(G)$
and $C(\gg ) = C'(\gg)$  for all 
$G\subset \SLV$. 
\ec 

\pf This is a direct consequence of Proposition \ref{i-iiProp}. 
\qed

In the following section we specialize our discussion to the case where
$G$ preserves a (possibly indefinite) scalar product or a symplectic form
on $V$.

\subsection{The orthogonal and symplectic cases} \label{ospSec} 
\bt Let $V$ be a pseudo-Euclidean vector space and $G \subset \SOV$ a Lie
subgroup. Then $A(G)$ is a point. In other words, the translation group $V$
is the only Abelian simply transitive affine group of type $G$. 
\et 

\pf By Theorem \ref{generalThm} we know that $A(G) \subset \gg^{(1)} \subset 
{\gso}(V)^{(1)}$. Now the theorem follows from the fact that 
${\gso}(V)^{(1)} = 0$. 
\qed

Let now $V$ be a symplectic vector space.  
The symplectic case is more interesting since the first prolongation of the 
symplectic Lie algebra is nontrivial. 
In fact, using the symplectic form $\o$ on $V$ we identify $V$ with $V^*$ 
via $v \mapsto \o v := \o (v, \cdot )$. This induces identifications   
of the Lie algebra $\gsp(V)$ with  the symmetric square $S^2V^*$ and of 
its  first prolongation with the symmetric cube $S^3V^*$, since  
\be \label{identification}  
\gsp(V)^{(1)} = V^* \ot \gsp(V) \cap
S^2V^*\ot V = V^* \ot S^2V^* \cap S^2V^* \ot V^* = S^3V^*\, .\ee
Explicitly, the identification of $S \in \gsp(V)^{(1)}$ with a 
totally symmetric trilinear form is given by $(X,Y,Z) \mapsto 
\o (S_XY,Z)$. Notice that $S^3V^*$ is also the vector space
of homogeneous cubic polynomials on $V$. The preceding discussion 
together with Corollary \ref{fundCor} establishes  therefore the following theorem. 

\bt The correspondence $H \mapsto S$ of Theorem \ref{generalThm} 
identifies $A(\SpV)$ with a quadratic cone $C(\gsp(V))$ 
in the vector space of cubic forms $S^3V^*$. 
\et  
  
In the following we identify of $S^3V^*$ with $S^3V$ by means
of $\o$.  Next we want to give an explicit construction of the solutions 
of the quadratic equations defining the affine variety 
$C(\gsp(V)) \subset S^3V$. 
We recall 
that the group $G$ acts by conjugation on $A(\gg)$. By 
Theorem \ref{generalThm}, under the identification $A(\gg)  = C(\gg )  
\hookrightarrow \gg^{(1)}$, this action is induced by  the natural linear
action on $\gg^{(1)}$. In the 
case of $G = \SpV$ this is the standard representation on $S^3V^* 
\cong S^3V$. We are also interested in the orbit space 
$\bar{A}(\gg) := A(\gg)/G$, 
which is the space of conjugacy classes of Abelian simply transitive affine
Lie algebras of type $G$. We have a natural inclusion 
$\bar{A}(\gsp(V) ) \subset S^3V/\SpV$.
 
\bt  \label{1stmainThm} Let $U \subset V$ be a Lagrangian subspace of a 
symplectic vector space $V$. Any $\SpV$-orbit  in
$C(\gsp(V))\subset S^3V$ intersects the subspace 
$S^3U \subset S^3V$. Moreover 
\[ C(\gsp(V)) = \bigcup_U S^3U \subset S^3V\, ,\] 
where the union
is over all Lagrangian subspaces $U \subset V$.
\et 

\pf Let  $S \in S^3 U$. Then the map $X \mapsto S_X$
takes values in $S^2U \subset  S^2V = \gsp(V)$. 
To see that $S \in C(\gsp(V))$ it is sufficient to observe that, since
$U$ is isotropic, $S^2U \subset  S^2V = \gsp(V)$ is an Abelian subalgebra
of linear operators.

Now let $S \in C(\gsp(V))$. We have to show that there exists a Lagrangian
subspace $U \subset V$ such that $S \in S^3U$. We prove first that the 
{\em support}  
\be \label{supportEqu} \S_S := {\rm span}  \{ S_XY \mid  X,Y \in V\}\ee 
of $S$ is isotropic. This implies $\S_S \subset U$ for some Lagrangian
subspace $U \in V$. It is sufficient to check that $S_X^2 = 0$ for all
$X\in V$. In fact, by polarization this implies $S_XS_Y = 0$ for all
$X,Y \in V$ and hence 
\[ \o (S_XY, S_ZW) = - \o (Y, S_XS_ZW) = 0\, .\]  
The claim that $S_X^2 = 0$ follows from the next computation. 
\[ \o (S_X^2Y,Z) = \o (S_YS_XX,Z) = -\o (S_XX,S_YZ) = 
-\o (S_XX,S_ZY) = \o (S_ZS_XX,Y) \] 
\[= \o (S_X^2Z,Y) = \o (Z,S_X^2Y) = 
-\o (S_X^2Y,Z) \, .\]
Now the theorem follows from the next lemma. 
\qed 

\bl \label{supportLemma}
Let $S\in S^3V$ and $\S_S \subset V_0$, where $V_0 \subset V$ is a 
subspace. Then $S \in S^3V_0$.
\el 

\pf Let $V_1 \subset V$ be a complement of $V_0$. We denote by 
$V_0^{\wedge}$ and $V_1^{\wedge}$ the annihilators of $V_0$ and $V_1$
in $V^*$, respectively, and consider them again as subspaces
of $V$. The decomposition 
\[ S^3V = \oplus_{p+q = 3} S^pV_0 \cdot S^qV_1 \, .\] 
gives rise to a decomposition $S = \sum_{p+q = 3} S^{p,q}$. For 
$X,Y\in V_0^{\wedge}$ we have $V_0 \ni S_XY = S^{0,3}_XY + 
S^{1,2}_XY$ and $S^{0,3}_XY \in V_1$. Hence, $S^{0,3} = 0$. 
Similarly, by considering $S_{V_1^{\wedge}}V_0^{\wedge}$ and  
$S_{V_1^{\wedge}}V_1^{\wedge}$ we can conclude that $S^{1,2} = 0$ and 
$S^{2,1} = 0$. This shows that $S = S^{3,0} \in S^3V_0$.  
\qed

For Lie subgroups of the  symplectic group we obtain: 
\bc Let $G \subset \SpV$ be a Lie subgroup. Then 
\[ C(\gg ) = \gg^{(1)} \cap
C(\gsp(V)) = \bigcup_U (\gg^{(1)}\cap S^3U) \subset S^3V\, ,\]
where the union
is over all Lagrangian subspaces $U \subset V$. 
\ec

Now we study more closely the cone $C(\gsp(V))$. Let us define 
\[ C(\gsp(V))_k := 
\{S\in  C(\gsp(V))\mid \dim \S_S = k\}\, . \]
 For any isotropic subspace $W\subset V$ we define 
\[ S^3W_{reg} := 
\{ S\in S^3W\mid \S_S = W\} \, .\]  
We denote the Grassmannian of isotropic subspaces $W\subset V$ of dimension $k$ by 
$G_k^0(V)$. Let ${\cal U}_k \ra G_k^0(V)$ be the universal (tautological) vector bundle. 

\bt \label{stratiThm} The cone $C(\gsp(V))$ has a stratification:
\[ C(\gsp(V)) = \bigcup_{k=0}^n C(\gsp(V))_k\, ,\]
where $2n = \dim V$. The $k$-th stratum 
\be \label{stratumEqu} 
C(\gsp(V))_k = \bigcup_{W \in G_k^0(V)} S^3W_{reg}  \ee
is a Zariski open subbundle of the vector bundle $S^3{\cal U}_k$. 
The action of $\SpV$ on $C(\gsp(V))$ preserves the stratification and we have the following
identification of orbit spaces:
\[ C(\gsp(V))_k/\SpV = S^3W_{reg}/{\rm GL}(W)\, ,\]
where $W  \subset V$ is a fixed $k$-dimensional isotropic subspace. 
\et

\pf It is clear that  
\[ C(\gsp(V)) = \bigcup_{k=0}^{2n} C(\gsp(V))_k\, .\]
It follows from Theorem \ref{1stmainThm} that the support $\S_S$ of any 
$S\in C(\gsp(V))$ is isotropic. This implies (\ref{stratumEqu}) and, 
in particular, that $C(\gsp(V))_k = \emptyset$ for all $k>n$.  
For the identification of orbit spaces it is sufficient to remark that two 
cubic forms $S$ and $S'$ with isotropic support $W$ are related by an 
element of $\SpV$ if and only if they are equivalent under ${\rm GL}(W)$.
In fact, since $W$ is isotropic, 
any element of ${\rm GL}(W)$ can be extended to an element of
$\SpV$, which preserves $W$. Conversely, any element of
$\SpV$ which maps $S$ to $S'$ has to preserve $W = \S_S = \S_{S'}$ and, hence,
induces an element of ${\rm GL}(W)$, which maps $S$ to $S'$.  
\qed 

\bc \label{translCor} Let $H \subset {\rm Aff}(\bR^{2n})$ be an Abelian 
simply transitive affine group of
symplectic type, $S\in S^3 V$ the cubic form defined by $H$.
Then the  subgroup of translations in $H$  is of dimension 
$2n-k$, where $k = \dim \S_S \le n$.  
\ec

\pf The subgroup of translations in $H$ is the
kernel of the map $X \mapsto S_X$ which coincides
with the kernel of $S$ as an element of $S^3 V^*$. Since
this is precisely the orthogonal complement of 
$\S_S$ with respect to $\omega$,  the dimension is $2n -k$.
Since $\S_S$ is isotropic, $k \leq n$. 
\qed 

  
\section{Flat special K\"ahler manifolds} \label{fskm}
\bd A (possibly indefinite) {\em K\"ahler manifold} is a (possibly indefinite) 
Riemannian manifold $(M,g)$ endowed with a parallel complex structure $J$, 
i.e.\ $DJ = 0$ for the Levi-Civita connection $D$ of $g$, so that $g$ is
Hermitian with respect to $J$.   
\ed
To any K\"ahler manifold $(M,g,J)$ we can canonically associate the parallel
symplectic form $\o = g(J \cdot, \cdot )$. It is called the 
{\em K\"ahler form}. 

\bd A {\em special K\"ahler manifold} is a K\"ahler manifold $(M,g,J)$ endowed
with a torsionfree and flat connection $\n$ such that
\begin{itemize}
\item[(i)] $\n \o = 0$ and 
\item[(ii)] $\n J$ is symmetric, i.e.\ $(\n_XJ)Y = (\n_YJ)X$ for all 
vector fields $X,Y$. 
\end{itemize}
A special K\"ahler manifold  $(M,g,J, \nabla)$ is called 
{\em flat\/} if the Levi-Civita
connection $D$ of $g$ is flat. 
\ed 

\noindent 
{\em Example:} Let $(M,g,J)$ be a {\em flat K\"ahler manifold}, 
i.e.\ a K\"ahler manifold for which the Levi-Civita connection 
$D$ is flat. Then $(M,g,J, D)$ is a flat special K\"ahler manifold. 
Special K\"ahler manifolds of this type will be called \label{trivial} 
{\em trivial special K\"ahler manifolds}.\\
 
The following characterization of trivial special K\"ahler manifolds
is easy to verify:

\bp Let $(M,g,J, \n )$ be a special K\"ahler manifold and $D$ its Levi-Civita
connection. Then the following three conditions are equivalent: 
\begin{itemize}
\item[(i)] $(M,g,J, \n )$ is a trivial special K\"ahler manifold, 
\item[(ii)]  $\n = D$ and  
\item[(iii)] $\n J = 0$. 
\end{itemize}
\ep 

Let us recall now some known results on special  K\"ahler manifolds
which relate to flatness. It was proven in \cite{BC} that any simply connected 
special K\"ahler manifold $(M,g,J,\n )$ of complex dimension $n$ 
can be canonically immersed into $\bR^{2n+1}$ as a parabolic affine
hypersphere with Blaschke metric $g$ and affine connection $\n$ (the 
notions of affine hypersphere, Blaschke metric and affine connection 
are discussed in detail in \cite{NS}). We recall the following classical 
theorem of 
Calabi and Pogorelov 
\cite{Ca}. 

\bt \label{CP}
If the Blaschke metric $g$ of a parabolic affine 
hypersphere $M$ of dimension $m$ is definite 
and complete, then 
$M$ is affinely congruent to the paraboloid $x^{m+1} 
= \sum_{i=1}^m (x^i)^2$ in $\bR^{m+1}$. 
In particular, $g$ is flat. 
\et
This implies: 

\bt \label{globalThm} 
Any special K\"ahler manifold $(M,g,J,\n )$ with a definite and
complete metric $g$ is trivial. 
\et 
 
\pf By \cite{BC}, the universal covering
is immersed as an affine hypersphere, which by Theorem \ref{CP} 
is a paraboloid. This implies that $\n = D$ is the Levi-Civita 
connection and hence that $(M,g,J,\n )$ is trivial. 
\qed   

Lu \cite{Lu} proved that any special K\"ahler manifold $(M,g,J,\n )$ 
with a definite and complete metric $g$ is flat.
Completeness is a very strong global assumption. 
Flatness on the other hand is a strong local assumption. 
If the metric is definite the following rigidity 
result holds:
\bt \label{rigidityThm} 
Any special K\"ahler manifold $(M,g,J,\n )$ with a definite and flat
metric is trivial. 
\et 

A proof will be given in the next section.
But let us remark that by the immersion theorem 
mentioned above Theorem \ref{rigidityThm} may also be read off 
from the classification of affine hyperspheres with flat and definite 
Blaschke metric \cite{VLS}.  
Without the assumption that the metric is definite 
Theorem  \ref{rigidityThm} 
does not hold. We will show that nontrivial flat indefinite 
special K\"ahler manifolds exist and give a local characterization
of such manifolds in section \ref{lc}. More specifically, 
we construct there examples of such manifolds which are geodesically complete 
with respect to both connections $D$ and $\n$. 

\subsection{The variety  $C_J(\gsp(V))$} \label{vc}
Let $(M,g,J,\n )$ be a flat special K\"ahler manifold.
We are interested in the local properties of $M$, and 
since $(M,g,J)$ is a flat K\"ahler-manifold we
may as well consider $M$ as ball in a standard
Hermitian vectors space $(V,J,g)$.  
The symplectic connection $\n$ may be expressed as
$\n  = D + S$, 
where  $S$ is a one-form on $M$ with values in $\gsp(V)$. We 
can think of it as a map $S: M \ra  V^* \ot S^2V^*$. 

\bp \label{sKProp} The tensor field $S$ satisfies the following
conditions
\begin{itemize}
\item[(i)] $S\in S^3V^* = \gsp(V)^{(1)}$, 
\item[(ii)] $[S_X,S_Y] = 0$  for all $X,Y\in V$,
\item[(iii)] $(D_XS)Y = (D_YS)X$ for all $X,Y \in V$ and 
\item[(iv)] $S_XJ = -J S_X$ for all $X\in V$.  
\end{itemize}
Conversely, any tensor field satisfying the conditions (i-iv) defines a flat 
special K\"ahler manifold $(M, g, J, \n := D+S)$. 
\ep 
Notice that (i-ii) say that $S$ has values in the cone 
$C(\gsp(V))$. \\

\pf The first condition is equivalent to the vanishing of the torsion of $\n$.
We claim that the Levi-Civita connection can be expressed as:  
\be D = \n - \frac{1}{2}J\n J \ee
The right-hand side is a torsionfree connection, by the symmetry of $\n J$.
Moreover it preserves $J$ and hence the metric $g = \o (\cdot , J\cdot )$: 
\begin{eqnarray*}(\n - \frac{1}{2}J\n J)_XJ &=& \n_XJ - \frac{1}{2}[J\n_XJ,J] 
\\ &=& \n_XJ - 
\frac{1}{2} J(\n_XJ)J +\frac{1}{2}  J^2\n_XJ = \frac{1}{2} (\n_XJ - J(\n_XJ)J)
= 0\, .\end{eqnarray*} 
For the last equation we have used that $J^2 = - {id}$ implies
\be \label{anticommEqu} J\n_XJ = -(\n_XJ)J \, .\ee
This proves the formula for the Levi-Civita connection, which we can restate
as: 
\be \label{SEqu} 2S_X = J\n_XJ \, .\ee
{}From (\ref{SEqu}) and (\ref{anticommEqu}) we conclude (iv).   
The remaining equations follow from the flatness of $\n$, as we will
show now. We compute the curvature $R^\n$ of $\n$: 
\begin{eqnarray*} R^\n (X,Y) &=& [\n_X,\n_Y] -\n_{[X,Y]} \\
&=& R^D (X,Y) + D_X(S_Y) - D_Y(S_X) +
[S_X,S_Y] -S_{[X,Y]}\\
&=& (D_XS)_Y - (D_YS)_X + [S_X,S_Y] \, .\end{eqnarray*} 
Here we have used that $D$ is flat, i.e.\ $R^D = 0$, and torsionfree. 
Now $[S_X,S_Y]$ commutes with $J$, whereas  $(D_XS)_Y - (D_YS)_X$
anticommutes with $J$, since $S_X$ and $S_Y$ anticommute with $J$, by (iv), 
and $DJ = 0$. This shows that $R^\n = 0$ implies (ii) and (iii).
 
Conversely, let $S$ be a tensor field satisfying (i-iv). Then 
$\n = D + S$ is a symplectic torsionfree (i) and by the formula above flat (ii-iii) connection.
It only remains to check the symmetry of $\n J$. For this we have 
to compute $\n J$ in terms of $S$ arriving again at (\ref{SEqu}):
\[ \n_XJ = D_XJ + [S_X,J] = [S_X,J] \stackrel{(iv)}{=} -2JS_X\, .\]
Now the symmetry of $S$ implies that of $\n J$. 
\qed  

The previous proposition shows that it is important to understand
the following closed conical subvariety of $\gsp(V)^{(1)} = S^3 V$: 
\be \label{CJEqu} C_J(\gsp(V)) := \{ S \in C(\gsp(V))\mid S_XJ + JS_X = 0 
\quad \mbox{for all} \quad X\in V\}\, .\ee 

\bl \label{isotropLemma} 
Let $S \in C_J(\gsp(V))$. Then $\S_S$ is $\o$-isotropic and 
$J$-invariant (and hence $g$-isotropic). 
\el 

\pf This follows immediately from Theorem \ref{1stmainThm}, the equation 
$JS_X = -S_XJ$ and the definition of the support $\S_S$, see equation 
(\ref{supportEqu}). 
\qed 

Now it is easy to prove Theorem \ref{rigidityThm}.\\

\pfof{Theorem \ref{rigidityThm}}
 The fundamental tensor field  
$S = \n -D$ associated to a flat special
K\"ahler manifold has values in $C_J(\gsp(V))$. We claim that
$C_J(\gsp(V)) = 0$ if the metric $g$ is definite. This implies 
the theorem. By the previous lemma, 
the support $\S_S$ of any $S \in C_J(\gsp(V))$ is $g$-isotropic and 
therefore $\S_S = 0$ if the metric is 
definite. 
\qed  

Next we want to describe $C_J(\gsp(V))$ for arbitrary signature of 
the metric. For this it it convenient to complexify $V$.  We obtain   
\be \label{cxifyEqu} V^{\bC} := V\otimes \bC = V^{1,0} \oplus V^{0,1} \, ,\ee
where $V^{1,0}$ and $V^{0,1}$ are the $\pm i$-eigenspaces of (the 
complex linear extension of) $J$. The symplectic form $\o$ extends to
a complex symplectic form (again denoted by $\o$) 
on $V^{\bC}$, for which the decomposition
(\ref{cxifyEqu}) is Lagrangian. We have a corresponding type decomposition 
of the third symmetric power: 
\[ S^3V^{\bC} = \oplus_{p+q =3} S^{p,q}V\, ,\quad S^{p,q}V := 
S^p V^{1,0}\cdot S^q V^{0,1}\, .\] 

\bl \label{typeLemma} 
Let $S \in S^3V$.  Then $S \in  C_J(\gsp(V))$ if and only if 
the following two conditions are satisfied
\begin{itemize}
\item[i)] the support $\S_S$ of $S$ is $\omega$-isotropic and $J$-invariant.
\item[ii)] $S = S_c + \overline{S}_c$, where
$S_c \in S^{3,0} V $ is a 
cubic form of type $(3,0)$ and $\overline{S}_c$ is the conjugated cubic 
form of type $(0,3)$. 
\end{itemize} 
Moreover, if $S \in  C_J(\gsp(V))$ then the support $\S_{S_c}$ of 
$S_c \in S^{3,0} V$ satisfies 
$ \S_{S_c} = \S_S^{1,0}$.     
\el  

\pf
If $S \in C_J(\gsp(V))$ condition i) is satisfied by Lemma \ref{isotropLemma}.
To show ii) we consider (see (\ref{identification})) $S$ as an element
of $S^3 ({V^{\bC}})^* = (S^3 V^{\bC})^*$. We put $V_{1,0}= (V^{0,1})^{\wedge}$
and $V_{0,1}= (V^{1,0})^{\wedge}$.  
We consider the induced dual decomposition 
$$S^3 ({V^{\bC}})^* = \oplus_{p+q =3} S_{p,q}V  ,\quad S_{p,q}V := 
S^p V_{1,0}\cdot S^q V_{0,1} \; .  
$$ 
The condition that 
$S_X$ and $J$ anticommute is expressed by the condition 
$S(J X, Y, Z) = S(X , JY, Z) = S(X, Y, JZ)$, for all $X, Y, Z \in V$.
This is equivalent to $$ S \in (S^{2,1} V  \oplus S^{1,2} V )^{\wedge}  
= S_{3,0} V \oplus  S_{0,3} V \, . $$ 
Under the identification via $\omega$
this amounts to $S \in  S^{3,0} V \oplus  S^{0,3} V$. 
Complex conjugation on $V^\bC$ extends to an  
antilinear involution $\rho$ on  $S^3 V^\bC$ so that $S^3V = (S^3W^\bC )^\rho$
is the fixed point set of $\rho$. 
Since $\rho$ interchanges $S^{3,0} V$ and  $S^{0,3} V$,  ii) holds. 
For the converse, we note that condition i) implies (see Lemma \ref{supportLemma})
that $S \in  C(\gsp(V))$. We just saw that condition ii) implies that 
$S_X$ and $J$ anticommute. Therefore i) and ii) imply $S \in  C_J(\gsp(V))$. 

We compute now the support  $\S_{S_c}$. 
$S_c \in S^{3,0} V = S^3 V^{1,0}$
implies that $\S_{S_c} \subset  V^{1,0}$, and correspondingly we get
$\S_{\overline{S}_c} \subset  V^{0,1}$. 
Since
$$ (\S_S)^{\bC} = (\S_S)^{1,0} \oplus (\S_S)^{0,1} = \S_{S_c} + \S_{\overline{S}_c}$$
we conclude that $\S_{S_c}=  \S_{S}^{1,0}$. 
\qed   

Let $W$ be a $J$-invariant isotropic subspace of $V$. Let us put
\[ S^{3,0}W^\rho := \{ S_c + \overline{S}_c\mid S_c \in S^{3,0}W \}
= (S^{3,0}W + S^{0,3}W)^\rho  \subset 
S^3W\, .\]    
The previous lemma implies the following theorem. 

\bt \label{3rdmainThm} 
\[ C_J(\gsp(V)) = \bigcup_U S^{3,0}U^\rho \subset S^3V\, ,\] 
where the union
is over all isotropic complex subspaces $U \subset V$ of maximal dimension. 
\et

Now we study more closely the cone $C_J(\gsp(V))$. Let us define 
\[ C_J(\gsp(V))_k := C_J(\gsp(V))\cap C(\gsp(V))_k \, .\]
We denote by $\bG_k(V)$ the Grassmannian of complex subspaces 
$W\subset V$ of (complex) dimension $k$ and by $\bG_k^0(V) \subset \bG_k(V)$
the real submanifold which consists of isotropic subspaces. 
Let $\bU_k \ra \bG_k^0(V)$ be the universal 
(tautological) vector bundle. It is the restriction of the holomorphic
universal bundle of the complex Grassmanian $\bG_k(V)$ to the 
submanifold $\bG_k^0(V) \subset \bG_k(V)$. 

\bt \label{JstratiThm} The cone $C_J(\gsp(V))$ has a stratification:
\[ C_J(\gsp(V)) = \bigcup_{k=0}^{N} C_J(\gsp(V))_k\, ,\]
where the number $N  \le n/2$ is the maximal complex dimension 
of an isotropic complex subspace of $V$. It is given by $N = \min (p,q)$ if 
$g$ has complex signature 
$(p,q)$, $p+q = n = \dim_{\bC}V$. The $k$-th stratum 
\be \label{stratumEqu2} 
C_J(\gsp(V))_k = \bigcup_{W \in \bG_k^0(V)} S^{3,0}W_{reg}^\rho 
\cong  \bigcup_{W \in \bG_k^0(V)} S^{3,0}W_{reg}\ee
is identified with a Zariski open subbundle of the vector 
bundle $S^3\bU_k$. 
The action of ${\rm Aut}(V,J,\o ) = {\rm Aut}(V,J,g) = {\rm U}(p,q)$ on 
$C(\gsp(V))$ preserves the stratification 
and we have the following
identification of orbit spaces:
\be \label{stratiEqu} C_J(\gsp(V))_k/{\rm U}(p,q) = 
S^{3,0}W_{reg}^\rho/{\rm GL}_{\bC}(W)
\, ,\ee
where $W  \subset V$ is a fixed complex $k$-dimensional isotropic subspace. 
\et 

\pf This follows from Theorem \ref{stratiThm} and Theorem \ref{3rdmainThm},  
taking under consideration the canonical identification  $S^{3,0}W^\rho \stackrel{\sim}{\ra} 
S^{3,0}W$, $S_c + \overline{S}_c \mapsto S_c$.
\qed 


\subsection{Local characterization of flat special special K\"ahler
manifolds} \label{lc}
Now we derive  classification results for indefinite special K\"ahler
manifolds. As before, let $(V,J,g)$ be a standard
(pseudo-) Hermitian vector space of complex signature $(p,q)$, $p+q = n$.  
 
\bt \label{4thmainThm} Let $f$ be a holomorphic function defined on an open 
subset $M \subset V^*$.
Assume that the pointwise support of the cubic tensor field 
$\partial^3f : M \ra  S^{3,0}V^*$ defined by the holomorphic  
third partial derivatives of $f$ is isotropic and put $S_f := \partial^3f + 
\overline{\partial^3f} \in S^3V^* \cong S^3V$. Then 
$M_f := (M,J,g, \n := D + S_f)$
is a flat special K\"ahler manifold and any flat special K\"ahler manifold
arises locally in this way. 
\et

\pf Let $f$ be a holomorphic function such that $\partial^3f$ has isotropic 
support. Then, by Lemma \ref{typeLemma},   the cubic tensor field 
$S_f := \partial^3f + \overline{\partial^3f} \in S^3V$ has also isotropic
support. This implies condition (ii) of Proposition \ref{sKProp}. Condition 
(iii) follows essentially from the construction of $S_f$ by means of partial  
derivatives, as we show next. We compute for constant (with respect to $D$) 
vector fields $X,Y,Z,W \in V$: 
\begin{eqnarray*} (D_XS)(Y,Z,W) &=& XS(Y,Z,W)\\
&=& \frac{1}{8}XY^{1,0}Z^{1,0}W^{1,0}f + \frac{1}{8}
XY^{0,1}Z^{0,1}W^{0,1}\overline{f}\\
&=& \frac{1}{16}X^{1,0}Y^{1,0}Z^{1,0}W^{1,0}f + \frac{1}{16}
X^{0,1}Y^{0,1}Z^{0,1}W^{0,1}\overline{f} 
 \, ,
\end{eqnarray*}  
where 
\[ X^{1,0} := \frac{1}{2}(X - iJ)\, ,\quad  X^{0,1} :=
\frac{1}{2}(X + iJ)\, .\]   
This shows that the tensor $DS$ is completely symmetric, which is the 
content of (iii) in Proposition \ref{sKProp}. Finally, (iv) 
is a direct consequence of Lemma \ref{typeLemma} and the definition
of $C_J(\gsp(V))$.    
So we have checked that $S_f$ satisfies the conditions (i-iv) of 
Proposition \ref{sKProp} and, hence, it defines a flat 
special K\"ahler manifold. 

Conversely, by Proposition \ref{sKProp}, any flat special K\"ahler manifold 
$(M,J,g, \n )$ is locally determined by the tensor field $S = \n -D : M \ra 
C_J(V)$, for which $DS$ is completely symmetric. This is the integrability
condition for the existence of a function $h$ such that $C = D^3h$. 
{}From Lemma \ref{typeLemma} it follows that
$S = S_c + \overline{S}_c \in S^{3,0}V + S^{0,3}V$ and the complete symmetry of
$DS$ implies the complete symmetry of $D^{1,0}S_c$ and $D^{0,1}\overline{S}_c$,
where $D = D^{1,0} +  D^{0,1}$ is the type decomposition of the (complexified)
flat torsionfree connection $D$. This shows the existence of a holomorphic
function $f$ such that $S_c = \partial^3 f$, where $\partial := D^{1,0}$. 
Thus, $S = S_f$. We know by Lemma \ref{isotropLemma} that the support
of $S$ is isotropic and, hence, by Lemma \ref{typeLemma}, also the support of 
$S_c$ is isotropic.
We remark that the  holomorphic function $f$ is related to 
the real function $h$ in the proof by $h = 8(f +\overline{f})$.  
\qed 

\noindent 
{\bf Remark 1} It follows from Proposition  \ref{sKProp}
that the tensor field $S = \n - D = S_f \in S^3V$ associated to a 
flat special K\"ahler manifold defines by means of $\o$ a commutative and 
{\em associative\/} multiplication $X\circ Y := S_XY$ on each tangent space. 
Moreover $S$ is potential in the sense that $S$ is defined by the third
derivatives of the function $h = 8(f +\overline{f})$ with respect to 
the flat torsionfree symplectic (and metric) connection $D$. 
This is similar to the type of structure one encounters in the theory of 
Frobenius manifolds. In that theory (see e.g.\ \cite{D}) 
one also has a cubic tensor field $S$ which is a section of $S^3T^*M$ 
and which is identified with a 
commutative and  associative multiplication by means of an isomorphism 
$T^*M \cong TM$. However, the isomorphism is given by a flat metric, 
whereas in our case we use the symplectic structure. In both cases 
$S$ is potential with respect to a flat torsionfree connection compatible 
with the identification $T^*M \cong TM$.  \\
 
\subsection{Examples with constant cubic form}
Let $(M,J,g,\n)$ be a special K\"ahler manifold. We say that it has {\em 
constant cubic form} if $DS = 0$, where $S = \n - D$.  
The special K\"ahler manifold $M_f$ in the above theorem has constant 
cubic form if and 
only if $f$ is a cubic polynomial. In that case we may assume that $f$ is 
a homogeneous polynomial, since lower order terms are annihilated by third derivatives. 
Notice that if $f$ is a homogeneous polynomial then the condition
on the third derivatives is satisfied if and only if $f$ considered
as an element of $S^3 V$ has isotropic support. 
We say that a special K\"ahler manifold $(M,J,g,\n )$ has a {\em trivial 
factor} if it is the product of two special K\"ahler manifolds
and one of the factors is trivial.
  
\bl 
The flat special K\"ahler manifold $(V,J,g,\n = D + S)$ defined by a cubic 
form $S\in C_J(\gsp(V))$ has a trivial factor if and only if 
$\dim \S_S < n = \frac{1}{2}\dim V$.    
\el 

\pf Suppose that $(V,J,g,\n = D + S)$ has a trivial factor. 
Then there exists an orthogonal and complex splitting $V = V_0 \oplus V_1$,
where $V_0 \subset \ker S$ is tangent to the trivial factor. 
Therefore $\S_S$ is contained in the nondegenerate subspace $V_1$.
This implies that the isotropic subspace $\S_S \subset V_1$ 
has dimension $\dim \S_S  \le \frac{1}{2}\dim V_1 < n$.
Conversely, if the isotropic space $W=\S_S$ has dimension $\dim W <n$
then we may choose an isotropic complement $W'$ of $W$ which
satisfies $\dim W' = \dim W$, so that $V_1 = W \oplus W'$ is
a nondegenerate subspace. Putting $V_0 = V_1^{\perp}$ we get
a complex orthogonal 
decomposition $V = V_0 \oplus V_1$, with $V_0 \subset \ker S$.
This shows that $V_0$ with
the induced structures is a trivial factor of $(V,J,g,\n = D + S)$. 
\qed

It follows from the Lemma and the remark above that any 
flat special K\"ahler manifold $M$ 
without a trivial factor and with constant cubic form is locally of the 
form $M_f$, 
where $f\in S^{0,3}U \subset S^{0,3}V \cong S^{3,0}V^*$ and 
$\Sigma_S = U \subset V$ is a complex Lagrangian
subspace. In particular the Hermitian vector space $(V,J,g)$ has
real signature $(2k,2k)$, and $\dim M = 4k$.   

\bt \label{5thmainThm} Let $(V,J,g)$ be a Hermitian vector space which
admits a complex Lagrangian subspace $U$, and let $f \in S^{0,3}U_{reg}
\subset S^{0,3}V \cong S^{3,0}V^*$. 
Then the manifold $M_f$ associated to the holomorphic cubic polynomial 
function $f: M = V \ra \bC$ is a  flat special 
K\"ahler manifold without trivial factor, with
constant cubic form and complete connections $D$ and $\n$. 
Conversely every flat special 
K\"ahler manifold without trivial factor and with
constant cubic form is locally of this form.  
Moreover, the correspondence $f \mapsto M_f$ defines a 
bijection between the orbits of the group ${\rm GL}_{\bC}(U)$ on 
$S^{0,3}U_{reg}$ and germs of flat  special 
K\"ahler manifolds with constant 
cubic form and without trivial factor up to equivalence.
\et    
\pf We already remarked that flat special K\"ahler manifolds
without trivial factor and constant cubic form are locally of the form 
$M_f$, for $f \in S^{0,3}U_{reg}$. 
Two germs of manifolds $M_f$ and $M_{f'}$ are isomorphic, 
if and only if there exists a germ of holomorphic isometry, i.e.\ an affine
transformation with linear part in ${\rm Aut}(V,J,g) = U(m,m)$ which    
maps the tensor $S_f$ to $S_{f'}$. In other words $M_f$ and $M_{f'}$ are 
isomorphic if and only if $f$ and $f'$ are equivalent under the (pseudo-) 
unitary group ${\rm U}(m,m)$.  {}From  Theorem \ref{JstratiThm} 
(\ref{stratiEqu}) we have the identification 
\[ \bigcup_{W\in \bG_m^0} S^3W_{reg}^\rho /{\rm U}(m,m) = 
C_J(\gsp(V))_m/{\rm U}(m,m) 
\cong S^{3,0}U_{reg}^\rho /{\rm GL}_{\bC}(U)
\cong S^{0,3}U_{reg}/{\rm GL}_{\bC}(U) \, ,\]
where the union is over all complex Lagrangian subspaces $W \subset V$. 
(They have complex dimension $N = m = n/2$.)  
This shows that the correspondence $S^{0,3}U_{reg}/{\rm GL}_{\bC}(U) \ni f
\mapsto M_f$ onto isomorphism classes of 
germs of flat  special K\"ahler manifolds with constant 
cubic form and without 
trivial factor is one-to-one. Note that $D$ is the canonical complete
connection on the vector space $M = V$. We explain now why $\n$ is complete.
Recall from section \ref{As} that by iii) of Proposition \ref{i-iiProp} 
the tensor $S = D -\n \in C(\gsp(V))$ defines 
a simply transitive affine action of $V$ on itself.
The orbit map $\Phi$ in $0$, $V \ni X \mapsto X + S_X X \in V$, 
of this action is therefore a diffeomorphism and it
is easy to very that $\nabla = \Phi^* D$. Hence, in particular $\nabla$
is complete.
\qed 
 
\noindent
{\bf Remark 2} The flat special K\"ahler manifolds 
$M_f = (V,J,g,\n = D + S_f)$ of Theorem \ref{5thmainThm} 
admit a simply transitive vector group of automorphisms. 
In fact, we can consider the constant tensor fields $J$, 
$g$ and $S_f$ as left-invariant complex structure, 
metric and connection on the vector group $V$. Every flat 
special K\"ahler manifold with a simply
transitive group of automorphisms is of this type.

\noindent
By the c-map, see \cite{ACD} and references therein, 
we can associate with the special K\"ahler manifold $M_f$ a 
flat hyper-K\"ahler manifold of signature $(4m,4m)$. 
This flat hyper-K\"ahler manifold then admits
a simply transitive group of automorphisms 
which is a semi-direct product of two vector 
groups of dimension $4m$.\\  

\noindent
{\bf Remark 3} The examples $M_f = (V,J,g,\n = D + S_f)$
are complete with respect to both connections $\nabla$
and $D$. It would be interesting to know if there
do exist further examples of flat special K\"ahler manifolds
(i.e. examples with a non-constant cubic form) 
which are complete with respect to both connections.

\end{document}